\documentclass{amsart}
\usepackage{amsmath,amssymb,latexsym,graphicx}

\newtheorem{theorem}{Theorem}

\begin{document}
\title{On the size of the resonant set for the products of $2\times 2$ matrices}
\author{Jeffrey Allen,
Benjamin Seeger,
Deborah Unger\\
}
\dedicatory{We would like to thank Professors Serguei Denissov and Alexander Kiselev for their constant support throughout the research and writing processes.  We are grateful to them for suggesting the topic and providing us with the opportunity to apply the concepts we have learned throughout our mathematical careers.  We thank them both for their excellent advice and encouragement.}
\maketitle

\begin{abstract}
    For $\theta \in [0, 2 \pi) $, consider the rotation matrix $R_\theta$ and
\[
h=\left(
\begin{array}{cc}
\lambda & 0\\
0 & 0
\end{array}\right), \quad \lambda>1.
\]
Let $W_n( \theta)$ denote the product of $m$ $R_{\theta}$'s and $n$ $h$'s with the condition $m\leq [\epsilon n]$ ($0 < \epsilon < 1$). We analyze the measure of the set of $\theta$ for which $\| W_n( \theta) \| \geq \lambda ^{\delta n}$ ($0 < \delta < 1$). This can be regarded as a model problem for the so-called Bochi-Fayad conjecture.
\end{abstract}

\section{Introduction}
In \cite{AR}, the following problem was considered. Take two matrices
\begin{equation}
H=\left(
\begin{array}{cc}
\lambda & 0\\
0 & \lambda^{-1}
\end{array}\right)\label{h}
\end{equation}
and
\[
R_{\theta} =
\begin{pmatrix}
    \cos (\theta) & - \sin (\theta) \\
    \sin (\theta) & \cos (\theta)
\end{pmatrix}, \quad \theta\in [0,2\pi) .
\]
Fix $\lambda>1$ and let $m, n \in \mathbb{N}$.
Define a word $W_n(\theta)$ to be
\[
    W_n(\theta) = H^{i_1} R_{\theta}^{j_1} \ldots H^{i_k} R_{\theta}^{j_k}
\]
where $i_1, \ldots , i_k, j_1, \ldots , j_k \in \mathbb{N} \cup \{0\}$ (where either or both $i_1$ and $j_k$ may be zero)
\[
    i_1 + \ldots + i_k = n \quad
    j_1 + \ldots + j_k = m
\]
and $k$ is arbitrary.
Assume that $m$ is much smaller than $n$ and take a ``generic'' angle $\theta$. It is not unreasonable to conjecture the geometric growth of  $\|W_n\|$ regardless of combinatorics of the word.
In \cite{AR}, the following theorem was proved:
\begin{theorem}
Assume that $0<\delta<1$ is fixed.
Then there is an $n$--independent set $\Omega$ such that $|\Omega|=2\pi$ and for any $\theta\in \Omega$ there is $\epsilon>0$ so that
\[
\min\limits_{W_n} \|W_n(\theta)\|>\lambda^{\delta n}
\]
provided
$m<\epsilon n (\ln n\ln\ln n)^{-1}$.
\end{theorem}

This theorem improved earlier results by Fayad and Krikorian \cite{FK}. The special case of the  Bochi-Fayad conjecture \cite{AR, FK} deals with the similar situation when $m<\epsilon n$ and $\epsilon$ is small. One might expect that $|\Omega|\to 2\pi$ as $\epsilon\to 0$ in this case. Proving it seems to be quite hard. We investigate a simpler case. In (\ref{h}), consider the matrix $H$  when $\lambda$ is large. Then $\lambda^{-1}\to 0$ as $\lambda\to \infty$ and one might wonder what happens if $\lambda ^ {-1}$ is dropped. We study the model case by taking
\[
h =
\begin{pmatrix}
    \lambda & 0 \\
    0       & 0
\end{pmatrix}\notin {\rm SL}(2,\mathbb{R})
\]
instead of $H$. It turns out, a very precise analysis can be performed. Section~3  provides some numerical evidence and comparison of the model case with  the real problem.
\section{The model problem}

In the previous setting, take

\[
h =
\begin{pmatrix}
    \lambda & 0 \\
    0       & 0
\end{pmatrix}
\]
instead of $H$ and $m \leq [\epsilon n]$, with $\epsilon\in (0,1)$ fixed.
Let
\[
W_n(\theta) =
\begin{pmatrix}
    a_n(\theta) & b_n(\theta) \\
    c_n(\theta) & d_n(\theta)
\end{pmatrix}
\]
and define a norm of $W_n(\theta)$ to be
\[
    \|W_n(\theta)\| = |a_n(\theta)| + |b_n(\theta)| + |c_n(\theta)| + |d_n(\theta)|.
\]
Set $$f_n(\theta) = \min_{W_n} \|W_n(\theta)\|.$$ Note that we can take the minimum because for an arbitrary $W_n$, there are $n$ $H$'s and $m$ $R_{\theta}$'s, so there are only a finite number of words for each $n$.
Finally, we fix $0 < \delta < 1$ and define the \textit{resonant set} $\mathcal{R}$ as thus: $\theta \in \mathcal{R}$ if there exists some $n$ such that $f_n(\theta) < \lambda ^{\delta n}$.
We claim that $|\mathcal{R}| < C \lambda ^{- \frac{1-\delta}{\epsilon}}$ where $C$ is some constant that can be explicitly computed and $|\mathcal{R}|$ denotes the Lebesgue measure of the set $\mathcal{R}$.

Given $n$ and $\theta$, there are four different types of words (note that here, none of $i_1, \ldots , i_{k+1}, j_1, \ldots , j_{k}$ are zero):
\begin{align}
    W_n(\theta) &= h^{i_1}R_{\theta}^{j_1} \cdots h^{i_k}R_{\theta}^{j_k} \label{E:HR} \\
    W_n(\theta) &= R_{\theta}^{j_1}h^{i_1} \cdots R_{\theta}^{j_k}h^{i_k} \label{E:RH} \\
    W_n(\theta) &= R_{\theta}^{j_1}h^{i_1} \cdots R_{\theta}^{j_{k-1}}h^{i_{k-1}}R_{\theta}^{j_k} \label{E:RR} \\
    W_n(\theta) &= h^{i_1}R_{\theta}^{j_1} \cdots h^{i_k}R_{\theta}^{j_k}h^{i_{k + 1}}. \label{E:HH}
\end{align}
In each word, there are precisely $k$ groups of rotation matrices of lengths $j_1, \ldots , j_k$. The only differences between the four types of words is which matrices ($h$ or $R_{\theta}$) they begin and end with.

From the word in (\ref{E:HR}), we obtain the following matrix:
\begin{align}
    W_n(\theta) &= h^{i_1} R_{\theta}^{j_1} \cdots h^{i_k} R_{\theta}^{j_k} \notag \\
                    &=
                    \begin{pmatrix}
                        \lambda^{i_1} & 0 \\
                        0             & 0
                    \end{pmatrix}
                    \begin{pmatrix}
                        \cos (j_1 \theta) & - \sin (j_1 \theta)\\
                        \sin (j_1 \theta) & \cos (j_1 \theta)
                    \end{pmatrix}
                        \cdots
                    \begin{pmatrix}
                        \lambda ^{i_k} & 0 \\
                        0              & 0
                    \end{pmatrix}
                    \begin{pmatrix}
                        \cos (j_k \theta) & - \sin (j_k \theta) \\
                        \sin(j_k \theta) & \cos (j_k \theta)
                    \end{pmatrix}
                    \notag \\
                    &=
                    \begin{pmatrix}
                        \lambda ^{i_1} \cos (j_1 \theta) & - \lambda ^{i_1} \sin(j_1 \theta) \\
                        0                                & 0
                    \end{pmatrix}
                        \cdots
                    \begin{pmatrix}
                        \lambda ^{i_k} \cos (j_k \theta) & - \lambda ^{i_k} \sin (j_k \theta) \\
                        0                                & 0
                    \end{pmatrix}
                    \notag \\
                    &=
                    \begin{pmatrix}
                        \lambda ^n \cos (j_1 \theta) \cdots \cos (j_k \theta) & \lambda ^n \cos (j_1 \theta) \cdots \cos (j_{k - 1} \theta) \sin (j_k \theta) \\
                        0                                                     & 0
                    \end{pmatrix}
                    .
                    \label{E:HRmatrix}
\end{align}
Likewise, from (\ref{E:RH}), we obtain
\begin{align}
    W_n(\theta) &= R_{\theta}^{j_1}h^{i_1} \cdots R_{\theta}^{j_k}h^{i_k} \notag \\
                &=
                    \begin{pmatrix}
                        \cos (j_1 \theta) & - \sin (j_1 \theta)\\
                        \sin (j_1 \theta) & \cos (j_1 \theta)
                    \end{pmatrix}
                    \begin{pmatrix}
                        \lambda^{i_1} & 0 \\
                        0             & 0
                    \end{pmatrix}
                        \cdots
                    \begin{pmatrix}
                        \cos (j_k \theta) & - \sin (j_k \theta) \\
                        \sin(j_k \theta) & \cos (j_k \theta)
                    \end{pmatrix}
                    \begin{pmatrix}
                        \lambda ^{i_k} & 0 \\
                        0              & 0
                    \end{pmatrix}
                    \notag \\
                 &=
                    \begin{pmatrix}
                        \lambda ^{i_1} \cos (j_1 \theta) & 0 \\
                        \lambda ^{i_1} \sin(j_1 \theta)  & 0
                    \end{pmatrix}
                        \cdots
                    \begin{pmatrix}
                        \lambda ^{i_k} \cos (j_k \theta) & 0 \\
                        \lambda ^{i_k} \sin(j_k \theta)  & 0
                    \end{pmatrix}
                    \notag \\
                 &=
                 \begin{pmatrix}
                    \lambda ^n \cos (j_1 \theta) \cdots \cos (j_k \theta)                  & 0 \\
                    \lambda ^n \sin (j_1 \theta) \cos(j_2 \theta) \cdots \cos (j_k \theta)  & 0
                 \end{pmatrix}
                 .
                 \label{E:RHmatrix}
\end{align}
Using the result in (\ref{E:RHmatrix}), the matrix from (\ref{E:RR}) is
\begin{align}
    W_n(\theta) &= R_{\theta}^{j_1}h^{i_1} \cdots R_{\theta}^{j_{k-2}}h^{i_{k-2}} R_{\theta}^{j_{k-1}}h^{i_{k-1}}R_{\theta}^{j_k} \notag \\
                 &=
                 \begin{pmatrix}
                    \lambda ^{i_1 + \ldots + i_{k-2}} \cos (j_1 \theta) \cdots \cos (j_{k-2} \theta)                  & 0 \\
                    \lambda ^{i_1 + \ldots + i_{k-2}} \sin (j_1 \theta) \cos(j_2 \theta) \cdots \cos (j_{k-2} \theta)  & 0
                 \end{pmatrix}
                 \notag \\
                 &
                 \quad \quad \quad \quad \quad \quad
                 \begin{pmatrix}
                    \cos (j_{k-1} \theta) & - \sin (j_{k-1} \theta) \\
                    \sin (j_{k-1} \theta) & cos(j_{k-1} \theta)
                 \end{pmatrix}
                 \begin{pmatrix}
                    \lambda ^ {i_{k-1}} & 0 \\
                    0                   & 0
                 \end{pmatrix}
                 \begin{pmatrix}
                    \cos (j_k \theta) & - \sin (j_k \theta) \\
                    \sin (j_k \theta) & \cos(j_k \theta)
                 \end{pmatrix}
                 \notag \\
                 &=
                 \begin{pmatrix}
                    \lambda ^{i_1 + \ldots + i_{k-2}} \cos (j_1 \theta) \cdots \cos (j_{k-2} \theta)                  & 0 \\
                    \lambda ^{i_1 + \ldots + i_{k-2}} \sin (j_1 \theta) \cos(j_2 \theta) \cdots \cos (j_{k-2} \theta)  & 0
                 \end{pmatrix}
                 \notag \\
                 &
                 \quad \quad \quad \quad \quad \quad
                 \begin{pmatrix}
                    \lambda ^{i_{k-1}} \cos (j_{k-1} \theta) \cos (j_k \theta)
                    &
                    - \lambda ^{i_{k-1}} \cos (j_{k-1} \theta) \sin (j_k \theta) \\
                    \lambda ^{i_{k-1}} \sin (j_{k-1} \theta) \cos (j_k \theta)
                    &
                    - \lambda ^{i_{k-1}} \sin (j_{k-1} \theta) \sin (j_k \theta)
                 \end{pmatrix}
                 \notag \\
                 &=
                 \begin{pmatrix}
                    \lambda ^n \cos (j_1 \theta) \cdots \cos (j_k \theta)
                    &
                    \lambda ^n \cos (j_1 \theta) \cdots \cos (j_{k - 1} \theta) \sin (j_k \theta)
                    \\
                    \lambda ^n \sin (j_1 \theta) \cos(j_2 \theta) \cdots \cos (j_k \theta)
                    &
                    \lambda ^n \sin (j_1 \theta) \cos(j_2 \theta) \cdots \cos(j_{k - 1} \theta) \sin (j_k \theta)
                 \end{pmatrix}
                 \label{E:RRmatrix}
\end{align}
and using the result in (\ref{E:HRmatrix}), the matrix from (\ref{E:HH}) is simply
\begin{align}
    W_n(\theta) &= h^{i_1}R_{\theta}^{j_1} \cdots h^{i_{k-1}}R_{\theta}^{j_{k-1}} h^{i_k}R_{\theta}^{j_k}h^{i_{k + 1}} \notag \\
                 &=
                 \begin{pmatrix}
                        \lambda ^{i_1 + \ldots + i_{k-1}} \cos (j_1 \theta) \cdots \cos (j_{k-1} \theta) & \lambda ^{i_1 + \ldots + i_{k-1}} \cos (j_1 \theta) \cdots \cos (j_{k - 2} \theta) \sin (j_{k-1} \theta) \\
                        0                                                     & 0
                 \end{pmatrix}
                 \notag \\
                 &
                 \quad \quad \quad \quad \quad \quad
                 \begin{pmatrix}
                    \lambda ^{i_k} & 0 \\
                    0              & 0
                 \end{pmatrix}
                 \begin{pmatrix}
                    \cos (j_k \theta) & - \sin (j_k \theta)\\
                    \sin (j_k \theta) & \cos (j_k \theta)
                 \end{pmatrix}
                 \begin{pmatrix}
                    \lambda ^{i_{k+1}} & 0 \\
                    0              & 0
                 \end{pmatrix}
                 \notag \\
                 &=
                 \begin{pmatrix}
                        \lambda ^{i_1 + \ldots + i_{k-1}} \cos (j_1 \theta) \cdots \cos (j_{k-1} \theta)
                        &
                        \lambda ^{i_1 + \ldots + i_{k-1}} \cos (j_1 \theta) \cdots \cos (j_{k - 2} \theta) \sin (j_{k-1} \theta) \\
                        0                                                     & 0
                 \end{pmatrix}
                 \notag \\
                 &
                 \quad \quad \quad \quad \quad \quad
                 \begin{pmatrix}
                        \lambda^{i_k + i_{k+1}} \cos(j_k \theta) & 0 \\
                        0                                        & 0
                 \end{pmatrix}
                 \notag \\
                 &=
                 \begin{pmatrix}
                   \lambda ^n \cos (j_1 \theta) \cdots \cos (j_k \theta) & 0 \\
                   0                                                     & 0
                 \end{pmatrix}
                 .
                 \label{E:HHmatrix}
\end{align}
Therefore we have
\[
    \|W_n(\theta)\| =
    \begin{cases}
        \lambda ^n | \cos(j_1 \theta) \cdots \cos(j_{k-1} \theta ) | (| \cos(j_k \theta)| + |\sin(j_k \theta)|), &\text{when $W_n$ satisfies (\ref{E:HR})}\\
        \lambda ^n |\cos(j_2 \theta) \cdots \cos(j_k \theta)| (|\cos(j_1 \theta)| + |\sin(j_1 \theta)|), &\text{when $W_n$ satisfies (\ref{E:RH})}\\
        \lambda ^n |\cos(j_2 \theta) \cdots \cos(j_{k-1} \theta)| ( |\cos(j_1 \theta) \cos(j_k \theta)| + |\cos(j_1 \theta) \sin(j_k \theta)|\\
        \, \, + |\sin(j_1 \theta) \cos(j_k \theta)| + |\sin(j_1 \theta) \sin(j_k \theta)| ), &\text{when $W_n$ satisfies (\ref{E:RR})}\\
        \lambda ^n |\cos(j_1 \theta) \cdots \cos(j_k \theta)| &\text{when $W_n$ satisfies (\ref{E:HH}).}
    \end{cases}
\]

{\bf Remark.}  This formula shows that $\min\limits_{W_n}\|W_n(\theta)\|$ is reached on the word of the type~(\ref{E:HHmatrix}).

\begin{theorem}{\label{T:resonant}}
    Let
    \[
    S_{\alpha} = \{ \, \theta \in [0,2 \pi ) \mid |\cos (\alpha \theta) | < \lambda ^{- \frac{(1 - \delta)\alpha}{\epsilon}-1 } \, \}
    \]
    and
    \[
    \tilde{S}_{\alpha} = \{ \, \theta \in [0,2 \pi ) \mid |\cos (\alpha \theta) | < \lambda ^{- \frac{(1 - \delta)\alpha}{\epsilon} } \, \}.
    \]
    Then the resonant set $\mathcal{R}$ satisfies
    $$\bigcup_{\alpha\in \mathbb{N}} S_\alpha\subseteq  \mathcal{R} \subseteq  \bigcup_{\alpha \in \mathbb{N}} \tilde{S}_{\alpha}.$$

\end{theorem}
\begin{proof}
    Suppose $\theta \in \bigcup_{\alpha \in \mathbb{N}} S_{\alpha}$. Then $\theta \in S_{\alpha}$ for some $\alpha \in \mathbb{N}$ and
    \[
        |\cos(\alpha \theta)| < \lambda ^{- \frac{(1-\delta) \alpha}{\epsilon}-1}.
    \]
    Let $n = [{\alpha/\epsilon}]+1$. Then, $n-1\leq \alpha/\epsilon<n$ and $\alpha<\epsilon n$. Consider the word $\omega _n(\theta) = h^{i_1} R_{\theta}^{\alpha} h^{i_2}$ where $i_1 + i_2 = n$. Since $m=\alpha$, we have $m\leq [\epsilon n]$. Then
    \[
        f_n(\theta) = \min_{W_n(\theta)} \|W_n(\theta)\| \leq \|\omega _n(\theta)\|
                    = \lambda ^n |\cos(\alpha \theta)| <  \lambda ^n \cdot \lambda ^{- \frac{(1-\delta) \alpha}{\epsilon}-1 } \leq  \lambda ^{\delta n}.
    \]
    Therefore $\theta \in \mathcal{R}$.

    Now suppose $\theta \notin \bigcup_{\alpha \in \mathbb{N}} \tilde{S}_{\alpha}$. Then for all $\alpha \in \mathbb{N}$,
    \[
        |\cos(\alpha \theta)| \geq \lambda ^{- \frac{(1-\delta) \alpha}{\epsilon}}.
    \]
    Choose arbitrary $n \in \mathbb{N}$. Then
    \begin{align*}
        f_n(\theta) &= \min_{W_n(\theta)} \|W_n(\theta)\| \\
                    &= \|\omega _n(\theta)\|
                        & &\text{(for some word $\omega _n(\theta)$) } \\
                    &= \lambda^n | \cos(j_1 \theta) \cdots \cos(j_k \theta) |
                        & &\text{(by the Remark above)} \\
                    &= \lambda^n | \cos (\alpha _1 \theta)^{m_1} \cdots \cos (\alpha _l \theta)^{m_l} |
    \end{align*}
    where $\alpha _1 < \ldots < \alpha_l$ and $m_1 \alpha _1 + \ldots + m_l \alpha _l = m \leq [\epsilon n]$. Then
    \begin{align*}
        f_n(\theta) &= \lambda^n |\cos (\alpha _1 \theta)^{m_1} \cdots \cos (\alpha _l \theta)^{m_l} |\\
                    &\geq \lambda^n \cdot \lambda^{- \frac{(1-\delta)(m_1 \alpha _1 + \ldots m_l \alpha_l)}{\epsilon}} \\
                    &= \lambda^n \cdot \lambda^{- \frac{m(1-\delta)}{\epsilon}} \geq  \lambda^{\delta n}
    \end{align*}
    and therefore $\theta \notin \mathcal{R}$.
\end{proof}

We claim that $\mathcal{R}$ is a dense open set. To show that $\mathcal{R}$ is open, we show that for each $n$, $f_n$ is continuous. For each $n$, $R_n = \{ \theta \in [0,2\pi) : f_n(\theta) < \lambda ^{\delta n} \} = f_n ^{-1} ((- \infty, \lambda ^{\delta n}))$, which is open as the pre-image of a continuous function of an open set. Note that $\mathcal{R} = \bigcup_{n=1}^{\infty} R_n$, a union of open sets, so $\mathcal{R}$ is open.

To show that $f_n$ is continuous, we note that $f_n$ is the minimum of a finite number of continuous functions (the norms of a finite number of words). Denote these functions as $F_1, F_2, \ldots , F_M$, $M \in \mathbb{N}$. Fix arbitrary $\theta \in [0, 2\pi)$, fix $\zeta > 0$, and let $\eta > 0$ be such that whenever $|\theta - \tilde{\theta}| < \eta$, $|F_k (\theta) - F_k (\tilde{\theta})| < \zeta$ for all $k = 1, \ldots, M$. Consider arbitrary $\tilde{\theta} \in (\theta - \eta, \theta + \eta)$. For some $i,j$, $f_n (\theta) = F_i (\theta)$ and $f_n (\tilde{\theta}) = F_j (\tilde{\theta})$. By the definition of $f_n$, $F_i (\theta) \leq F_j (\theta)$ and $F_j (\tilde{\theta}) \leq F_i (\tilde{\theta})$. Notice that if $F_i(\theta) = F_j(\tilde{\theta})$, then $|f_n (\theta) - f_n (\tilde{\theta})| = 0 < \zeta$ and we are done. Suppose that $F_i(\theta) > F_j(\tilde{\theta})$. Then $|f_n (\theta) - f_n (\tilde{\theta})| = F_i(\theta) - F_j(\tilde{\theta}) \leq F_j(\theta) - F_j(\tilde{\theta}) < \zeta$. Otherwise, if $F_i(\theta) < F_j(\tilde{\theta})$, then $|f_n (\theta) - f_n (\tilde{\theta})| = F_j(\tilde{\theta}) - F_i(\theta) \leq F_i(\tilde{\theta}) - F_i(\theta) < \zeta$.

To see that $\mathcal{R}$ is dense, let $I$ be any open interval in $[0, 2 \pi)$. The collection of points $R_{\alpha} = \{\frac{\pi}{2 \alpha} + \frac{\pi}{\alpha} k : k \in \{1, \ldots , 2 \alpha - 1 \} \}$ is in $S_{\alpha}$; indeed, for any $\phi \in R_{\alpha}$, $\cos(\alpha \phi) = 0 < \lambda ^{- \frac{(1 - \delta)\alpha}{\epsilon}-1 }$. If we choose $\alpha > \frac{|I|}{\pi}$, then there must be some element $\phi$ of $R_{\alpha}$ in $I$. Since $\phi \in \bigcup_{\alpha\in \mathbb{N}} S_\alpha\subseteq  \mathcal{R} \subseteq  \bigcup_{\alpha \in \mathbb{N}} \tilde{S}_{\alpha}$, we see that every open interval in $[0, 2 \pi)$ contains a point in $\mathcal{R}$.

Now we are ready to estimate the size of $\mathcal{R}$. Consider $\tilde{S}_{\alpha}$ for arbitrary $\alpha \in \mathbb{N}$. The measure of this set is
\begin{align*}
    |\tilde{S}_{\alpha}| &= 4 \alpha (\frac{\pi}{2 \alpha} - \frac{1}{\alpha} \cos ^{-1} ( \lambda ^ {- \frac{(1-\delta)\alpha}{\epsilon}}))\\
             &= 2 \pi - 4 \cos ^{-1} ( \lambda ^ {-\frac{(1-\delta) \alpha}{\epsilon}}) \\
             &\approx 2\pi - 2\pi + 4 \lambda ^ {- \frac{(1-\delta) \alpha}{\epsilon}} \\
             &= 4 \lambda ^ {- \frac{(1-\delta) \alpha}{\epsilon}}.
\end{align*}

Then our estimate for the size of $\mathcal{R}$ is
\begin{align*}
    |\mathcal{R}| &\leq |\bigcup_{\alpha \in \mathbb{N}} \tilde{S}_{\alpha}| \leq \sum_{\alpha=1}^{\infty} |\tilde{S}_{\alpha}| \approx 4\sum_{\alpha=1}^{\infty}  \lambda ^ {- \frac{(1-\delta)\alpha}{\epsilon}} \\
           &= \frac{4 \lambda^{- \frac{1-\delta}{\epsilon} } }{1 - \lambda^{- \frac{1-\delta}{\epsilon}}} \approx 4 \lambda^{- \frac{1-\delta}{\epsilon}}
\end{align*}
as $\epsilon\sim 0$ and $\lambda$, $\delta$ are fixed.
\section{Some numerical evidence}
We provide some numerical and graphical evidence of what was proved. We see how the graphs of the model case and the real case compare for fixed $n$ and $m$. In addition, it is shown graphically that changing the multiplicities of $H$ affects the word's norm, but in the model case the word's norm is invariant under these changes.  The graphs in this section were plotted with Maple 14.

Based on the similarities between the pictures of the real case and the model case, we conjecture that the resonant set in the model case is, in some sense, the limiting set of the the resonant set in the real case as $\lambda$ grows large, since the $\lambda^{-1}$ term goes to $0$ as $\lambda$ goes to infinity. Of course, since here we take $\lambda$ relatively small ($\lambda = 2$) for graphing convenience, this is a rough conjecture; in fact, proving it seems to be rather difficult.

\begin{figure}[h]
\begin{center}
\includegraphics[height=5cm,width=10cm]{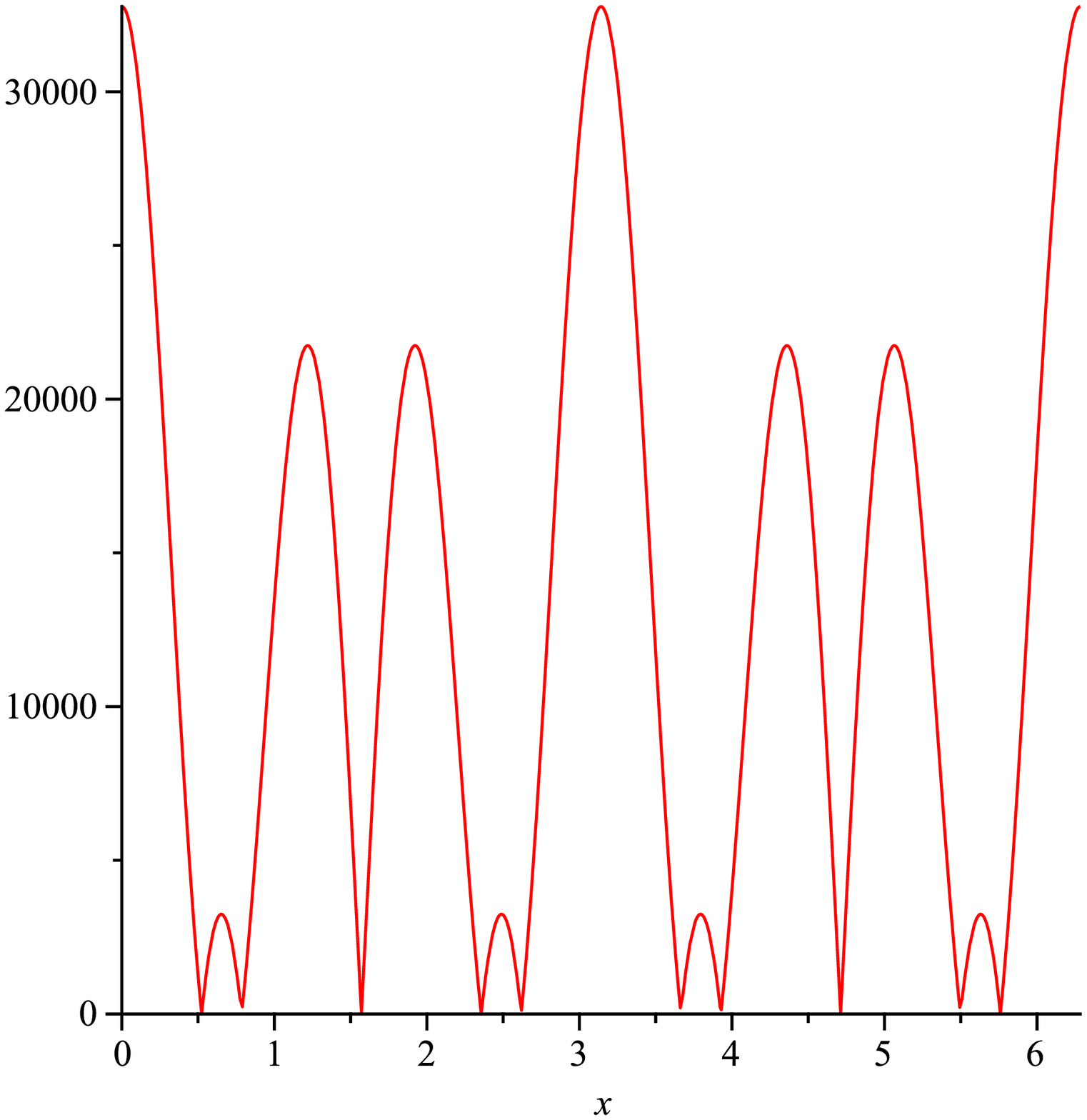}
\end{center}
\caption{}\label{f1}
\end{figure}
Figure \ref{f1} shows the graph of $||h^{i_1} R_{\theta}^2 h^{i_2} R_{\theta}^3 h^{i_3}||$ where $i_1,i_2,i_3 \in \mathbb{N}$ and $i_1 + i_2 + i_3 = 15$. Recall that $h$ is the matrix we use in the model case, where $\lambda ^{-1}$ is replaced by 0. With these combinatorics, varying $i_1,i_2,i_3$ does not change the graph, as long as their sum is 15.  Specifically, Figure \ref{f1} is the model case of $n = 15$, $m = 5$, and $\lambda = 2$.

\begin{figure}[h]
\begin{center}
\includegraphics[height=5cm,width=10cm]{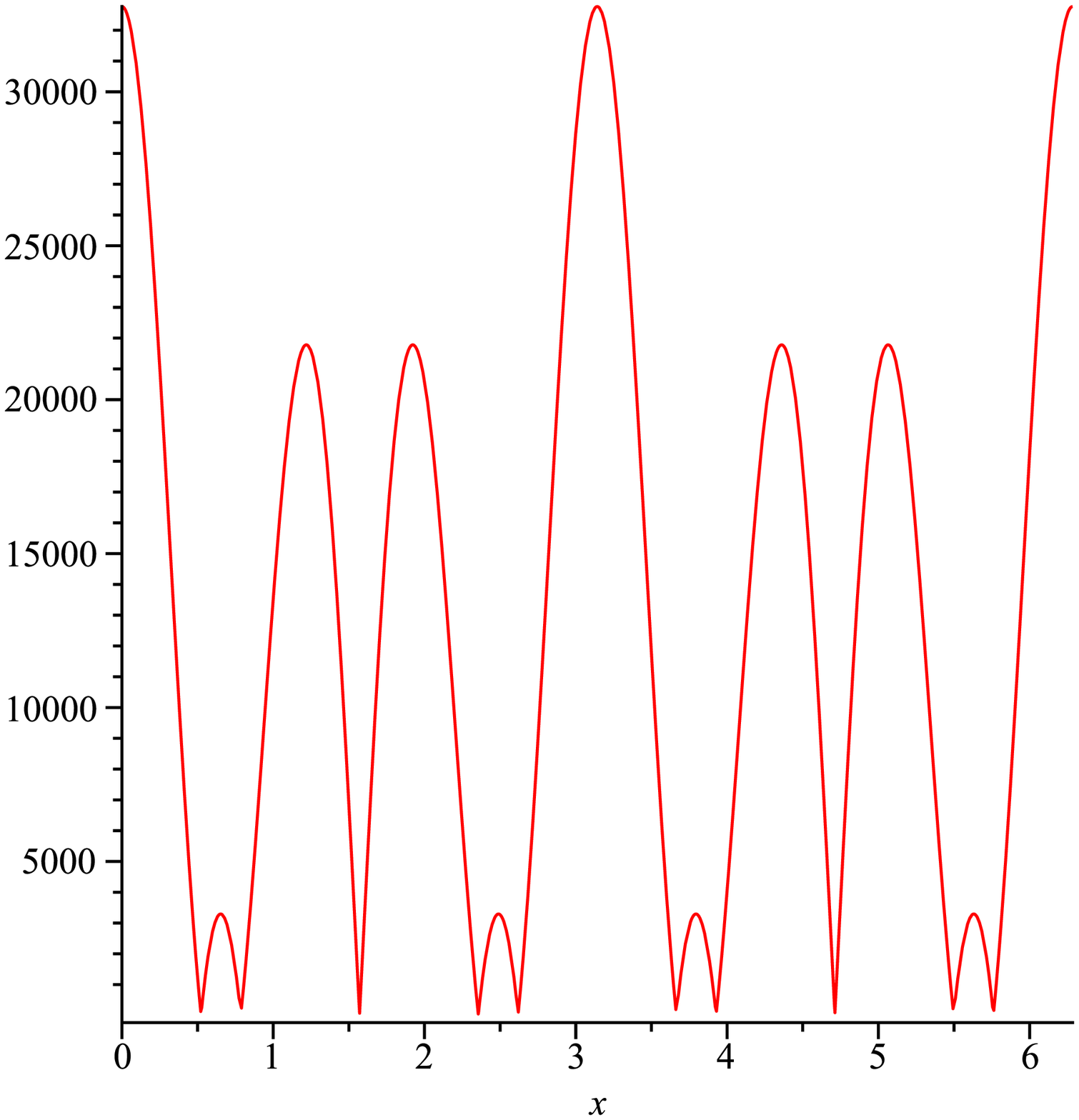}
\end{center}
\caption{}\label{f2}
\end{figure}
Figure \ref{f2} shows the case where we replace $h$ with $H$ and set $i_1=i_2=i_3=5$.  Explicitly, we are graphing $|| H^5 R_{\theta}^2 H^5 R_{\theta}^3 H^5 ||$.

\begin{figure}[h]
\begin{center}
\includegraphics[height=5cm,width=10cm]{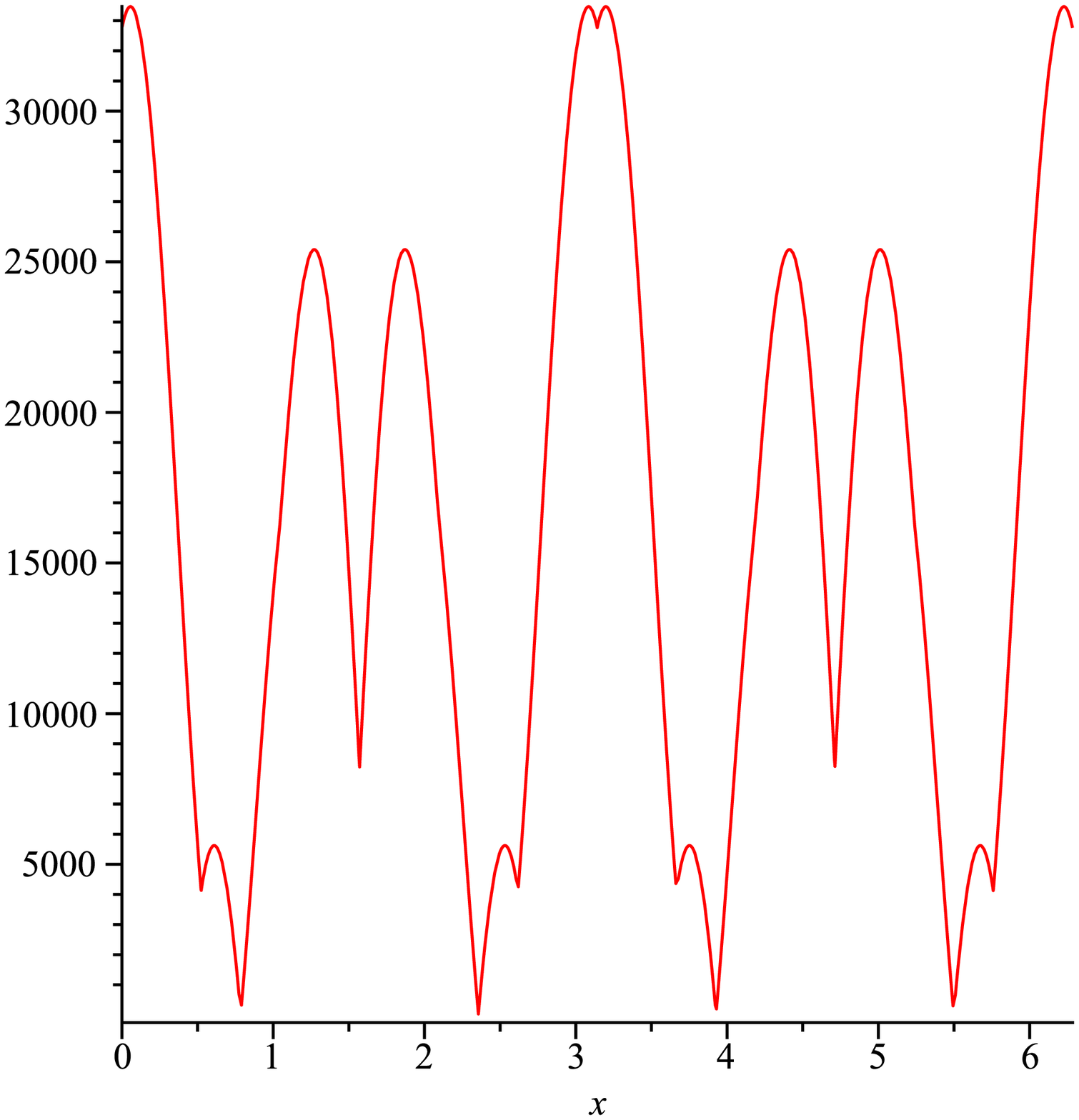}
\end{center}
\caption{}\label{f3}
\end{figure}
Figure \ref{f3} shows the graph of $|| H^5 R_{\theta}^2 H^9 R_{\theta}^3 H^1 ||$, where every variable is the same as above, $n = 15$, $m = 5$, and $\lambda = 2$, but the $H$ matrices are multiplied in different orders with respect to the rotation matrices. We still have $i_1 + i_2 + i_3 = n = 15$, however, $i_2 = 9$ and $i_3 = 1$.  Changing the order of multiplication in the word with the same number of $H$ matrices changes the norm of the word.

By comparing Figure \ref{f2} and Figure \ref{f3}, we observe that a greater disparity between the multiplicities of $H$ (the $i_k$'s) is correlated with a smaller resonant set (the set of points $\theta$ between $0$ and $2 \pi$ such that the norm of the word is within a certain distance of zero).  The slope of the word's norm is steeper in Figure \ref{f3} than it is in Figure \ref{f2} and the peaks in Figure \ref{f3} are associated with larger values of the word's norm than in the case depicted by Figure \ref{f2}. Both conditions lead to fewer points $\theta$ that are mapped to a norm of the word that is close to zero.

\begin{figure}[h]
\begin{center}
\includegraphics[height=5cm,width=10cm]{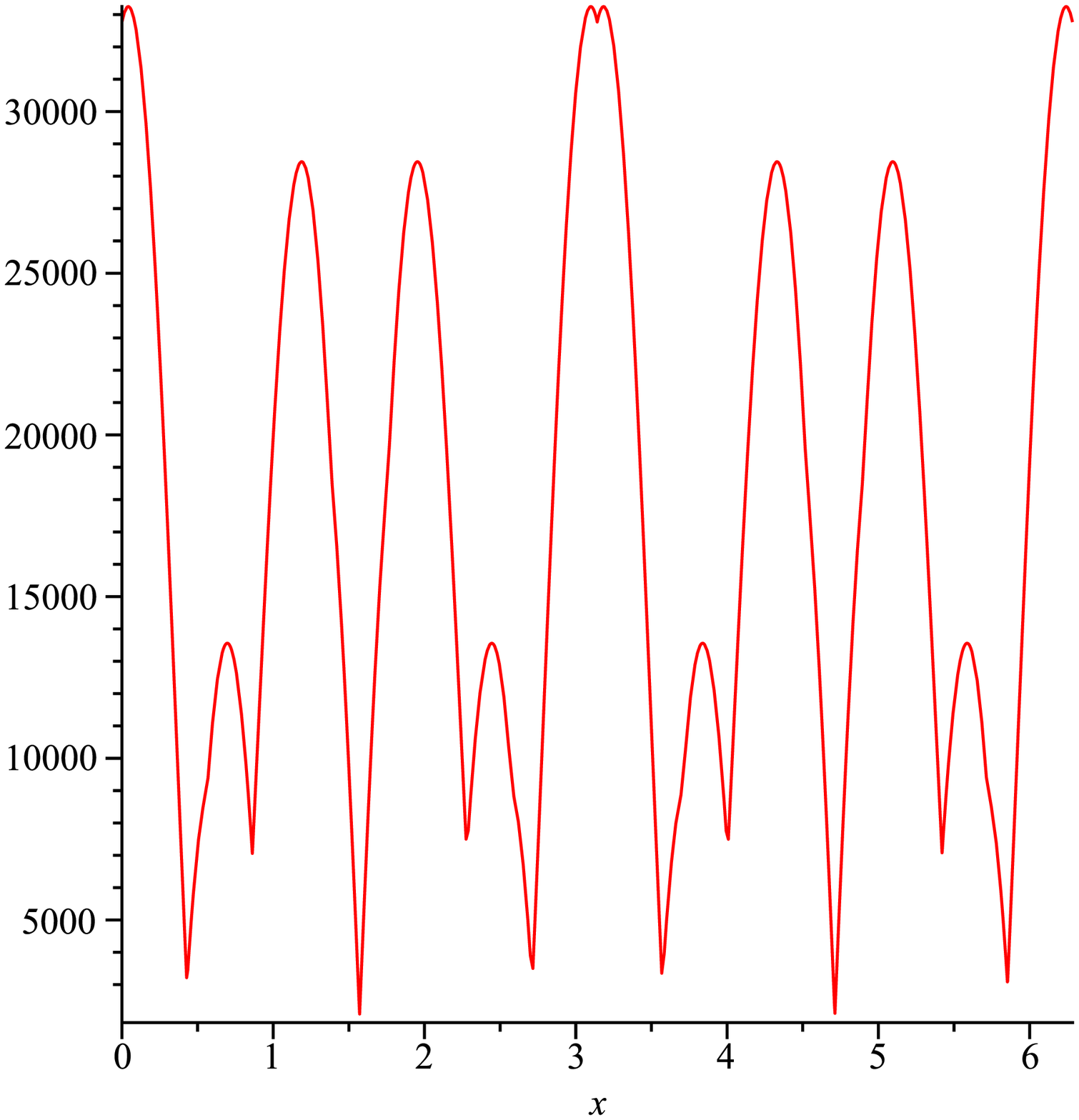}
\end{center}
\caption{}\label{f4}
\end{figure}
Figure \ref{f4} shows the graph of $|| H^1 R_{\theta}^2 H^1 R_{\theta}^3 H^{13} ||$.  Comparing this graph with Figure \ref{f3} provides further evidence that a greater disparity between the multiplicities of $H$ results in a smaller resonant set.


\begin{figure}[h]
\begin{center}
\includegraphics[height=5cm,width=10cm]{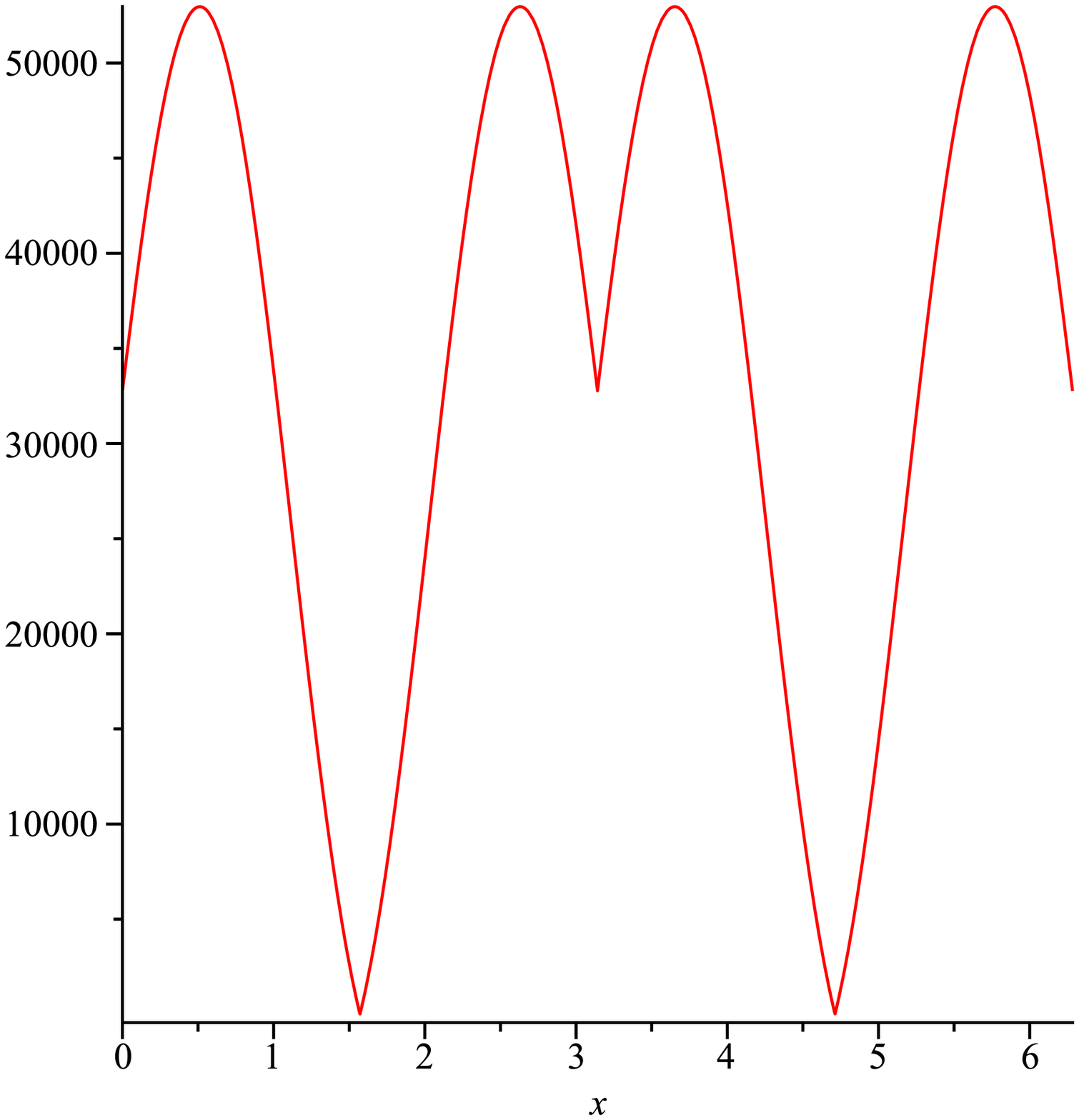}
\end{center}
\caption{}\label{f5}
\end{figure}
\begin{figure}[h]
\begin{center}
\includegraphics[height=5cm,width=10cm]{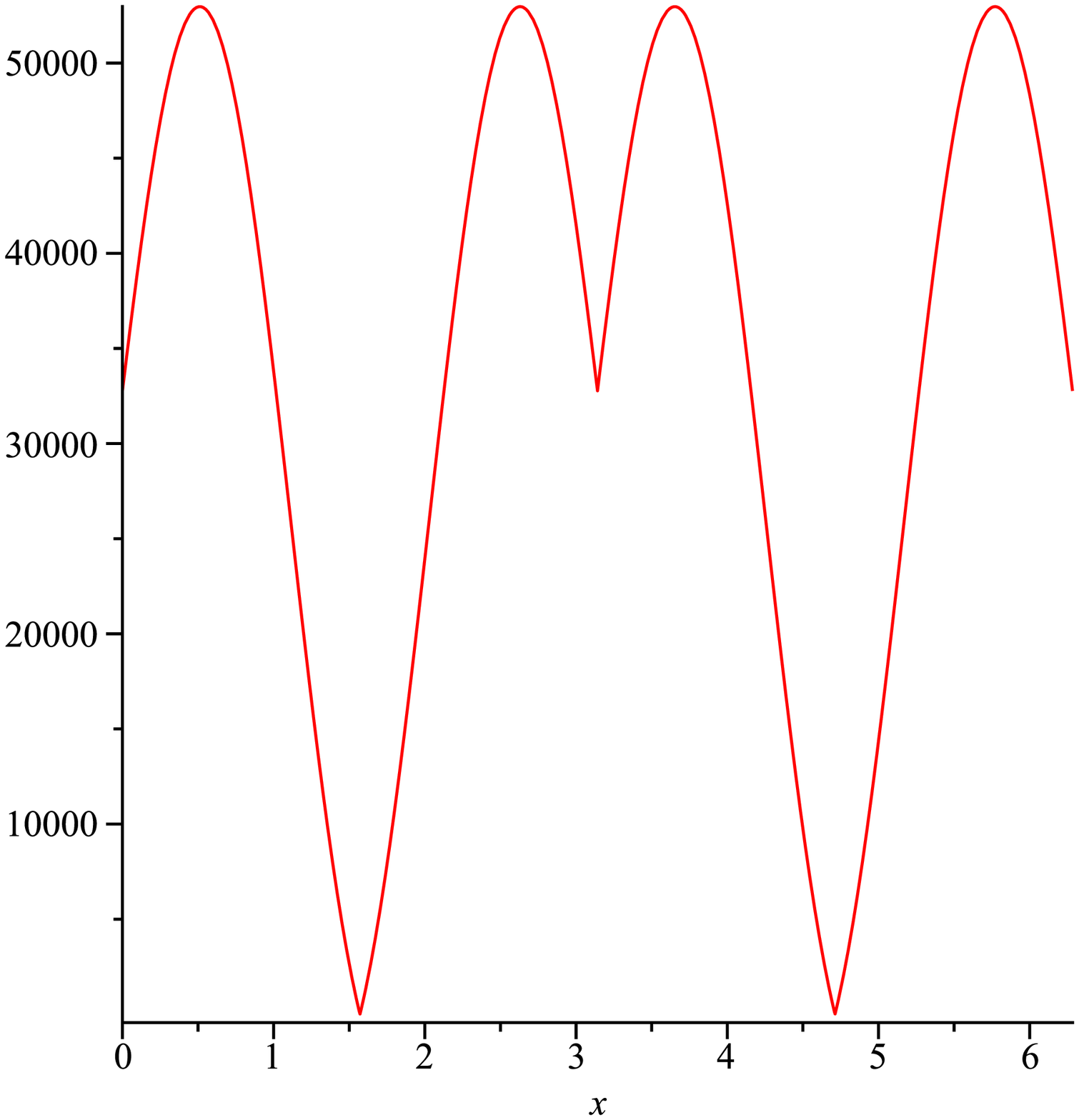}
\end{center}
\caption{}\label{f6}
\end{figure}

To further justify our use of the model case where the fourth entry in $H$ is 0 and $\lambda = 2$, notice that the graph in Figure \ref{f5} of the model case, $|| R_{\theta} h^{i_1} R_{\theta} h^{i_2} R_{\theta} ||$ where $i_1 + i_2 = 15$, and the graph in Figure \ref{f6} of the real case, $|| R_{\theta} H^7 R_{\theta} H^8 R_{\theta} ||$, appear to be nearly identical. As proved above, the word is independent of the order of multiplication of the $h$'s so long as $k = 2$ and $n = 15$.

Furthermore, we notice that Figure \ref{f5} is not comparable to Figure \ref{f1} although both show the model case, but with different combinatorics on the word.  Figure \ref{f1} represents an example of equation 5, whereas Figure \ref{f5} represents an example of equation 4.  Both graphs still have a small resonant set.

\begin{figure}[h]
\begin{center}
\includegraphics[height=5cm,width=10cm]{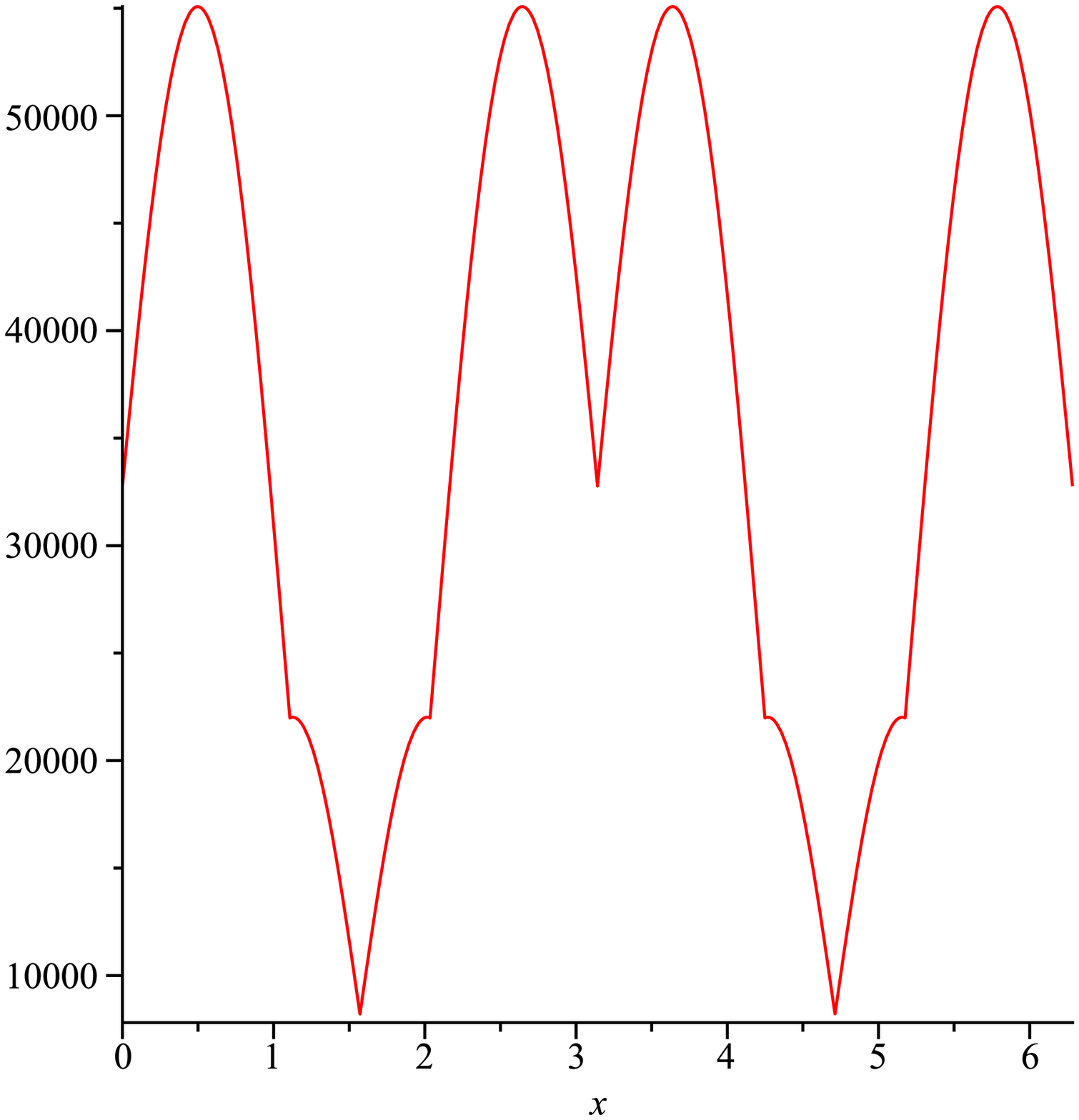}
\end{center}
\caption{}\label{f7}
\end{figure}
Figure \ref{f7} shows the graph of $|| R_{\theta} H^{14} R_{\theta} H R_{\theta} ||$.  By comparing Figure \ref{f6} and Figure \ref{f7}, again we see that a greater disparity between the multiplicities of $H$ results in a smaller resonant set.

\end{document}